\documentstyle[11pt]{article}

\input amssym.def
\input amssym

\def \a{{\frak a}}
\def \A{{\Bbb A}}
\def \Ad{{\rm Ad}}
\def \bs{{\backslash}}
\def \CG{{\cal G}}
\def \CH{{\cal H}}
\def \C{{\Bbb C}}
\def \CO{{\cal O}}
\def \CR{{\cal R}}

\def \g{{\frak g}}
\def \ga{\gamma}
\def \Ga{{\Gamma}}

\def \Hom{{\rm Hom}}
\def \ind{{\hspace{2pt}\rm ind\hspace{1pt}}}
\def \k{{\frak k}}
\def \la{{\lambda}}

\def \Lie{{\rm Lie}}

\def \p{{\frak p}}
\def \ph{\varphi}
\def \prf{{\bf Proof: }}

\def \qed{\hfill $\Box$

$ $

}
\def \Q{{\Bbb Q}}
\def \ra{\rightarrow}
\def \rank{{\hspace{2pt}{\rm rank}\hspace{1pt}}}
\def \R{{\Bbb R}}
\def \Res{{\rm Res}}
\def \t{{\frak t}}
\def \tr{{\hspace{2pt}\rm tr\hspace{1pt}}}
\def \vol{{\hspace{1pt}\rm vol\hspace{1pt}}}

\begin{document}

\title{On the index of Dirac operators on arithmetic quotients}

\author{Anton Deitmar\\ $ $\\ {\small Math. Inst. d. Univ., INF 288, 69126
Heidelberg, Germany}}
\date{}
\maketitle

\pagestyle{myheadings}
\markright{ON THE INDEX...}

\begin{scriptsize}

{\bf ABSTRACT:}
Using the Arthur-Selberg trace formula we express the index of a Dirac operator
on an arithmetic quotient manifold as the integral over the index form plus a
sum of orbital integrals.
For the Euler operator these orbital integrals are shown to vanish for products
of certain rank one spaces.
In this case the index theorem looks exactly as in the compact case.

\end{scriptsize}

\tableofcontents

\newcommand{\rez}[1]{\frac{1}{#1}}
\newcommand{\der}[1]{\frac{\partial}{\partial #1}}
\newcommand{\binom}[2]{\left( \begin{array}{c}#1\\#2\end{array}\right)}

\newcounter{lemma}
\newcounter{corollary}
\newcounter{proposition}
\newcounter{theorem}

\newtheorem{conjecture}{\stepcounter{lemma} \stepcounter{corollary}
	\stepcounter{proposition}\stepcounter{theorem}Conjecture}[section]
\newtheorem{lemma}{\stepcounter{conjecture}\stepcounter{corollary}
	\stepcounter{proposition}\stepcounter{theorem}Lemma}[section]
\newtheorem{corollary}{\stepcounter{conjecture}\stepcounter{lemma}
	\stepcounter{proposition}\stepcounter{theorem}Corollary}[section]
\newtheorem{proposition}{\stepcounter{conjecture}\stepcounter{lemma}
	\stepcounter{corollary}\stepcounter{theorem}Proposition}[section]
\newtheorem{theorem}{\stepcounter{conjecture} \stepcounter{lemma}
	\stepcounter{corollary}	\stepcounter{proposition}Theorem}[section]

$$ $$

\begin{center} {\bf Introduction} \end{center}

Index theorems for noncompact but finite volume locally symmetric spaces $Y$
are usually of the form
$$
\ind (D) = \int_Y \omega + "{\rm error\ terms}",
$$
where $\omega$ is the local index form of the elliptic differential operator
$D$.
The error terms are associated to the cusps.

F. Hirzebruch \cite{Hirz} used resolution of cusp-singularities to show that in
the case of Hilbert modular surfaces and their signature operators the error
terms can be given as special values of $L$-functions.
Using the Selberg trace formula this was extended to higher dimensional Hilbert
modular varieties by W. M\"uller \cite{mull}.
Also by means of the Selberg trace formula D. Barbasch and H. Moscovici
\cite{BM} showed index theorems for real rank one spaces.
A detailed analysis of the geometry of arithmetic quotient varieties led M.
Stern \cite{stern1}, \cite{stern2} to index theorems of the above type.
By means of the adelic trace formula J.P. Labesse \cite{Lab-sig} investigated
the index of the signature operator focussing on the representation theoretic
aspects.

In this paper we will use Arthur's formula \cite{Art-inv2} to give a geometric
index formula of the above type for Dirac operators.
The advantage of this index formula is that in some cases like products of real
hyperbolic spaces the vanishing of the error term can be read off.
The most important example for this is the Euler operator $D^e =d+d^*$ as
operator from even to odd forms.
In the compact case one has
$$
\ind D^e = \int_Y \omega = \chi (Y),
$$
where $\omega$ is the Euler form and $\chi (Y)$ the topological Euler
characteristic.
The first equality is the Atiyah-Singer index theorem, the second the
Gauss-Bonnet equality.
In the case of arithmetic quotients G. Harder \cite{Hard} has shown that the
second equality also holds in the noncompact case.
One might wonder whether the first also carries over to noncompact arithmetic
quotients.
In this paper we show that the first equality holds in case of products of rank
one spaces modulo arithmetic groups coming from totally real number fields with
at lest two real embeddings.
This contains the case of Hilbert modular varieties.

\section{Homogeneous Dirac operators and their index kernels}
An {\bf arithmetic quotient} is a quotient manifold $X_\Ga =\Ga \bs X$ of a
globally symmetric space $X$ by a torsion-free arithmetic group $\Ga$.
The space $X$ is assumed not to have compact or euclidean factors.
Under these circumstances the space $X$ can be written as homogeneous space
$X=G/K$, where $G$ is a semisimple real reductive group acting transitively by
isometries and $K$ is a maximal compact subgroup of $G$.
The group $\Ga$ is an arithmetic subgroup of $G$.

To be able to use adelic methods we will further assume that $\Ga$ is a {\bf
congruence subgroup}, i.e. there is a semisimple linear algebraic group $\CG$
over $\Q$ with $G=\CG(\R)$, a compact open subgroup $K_\Ga$ of $\CG(\A_{fin})$,
where $\A_{fin}$ is the ring of finite adeles over $\Q$ such that $\Ga =
\CG(\Q)\cap K_\Ga$.

In order to have strong approximation available we will also assume the group
$\CG$ to be simply connected.

Since the indices of homogeneous Dirac operators are known to vanish otherwise
we will assume
$$
{\rm rank}\ G = {\rm rank}\ K.
$$

On the group $\CG(\A)$ we have a unique Haar measure given by a rational top
differential form.
This measure is also called the {\bf Tamagawa measure}.
We will distribute the Tamagawa measure to the factors $\CG(\R)$ and
$\CG(\A_{fin})$ in a way that on $G=\CG(\R)$ we have the {\bf Euler-Poincar\'e
measure} given by
$$
\vol (\Ga \bs G) = (-1)^{\frac{\dim X}{2}}\chi (X_\Ga),
$$
where $\chi$ denotes the Euler-Poincar\'e characteristic. It suffices to insist
this formula to hold for cocompact torsion free lattices $\Ga$ but for
arithmetic ones which are not cocompact it holds as well (\cite{Hard}).

Let $\g_0 =\k_0 \oplus \p_0$ denote the polar decomposition of the real Lie
algebra $\g_0$ of $G$ where $k_0 := \Lie_\R K$ and $\p_0$ is the orthogonal
complement of $\k_0$ with respect to the Killing form $B$ of $\g_0$.
Let $\g = \k \oplus \p$ denote the complexified version.
Since $G$ and $K$ have the same rank there is a Cartan subgroup $T$ of $G$
which is contained in $K$.
Let $\t$ denote the complex Lie algebra of $T$.
Choose an ordering on the root system $\Phi (\g ,\t)$.
Since $\p$ is stable under $\Ad(T)$ it follows that this choice induces a
decomposition $\p = \p_+ \oplus \p_-$ according to positive and negative root
spaces.
This decomposition is a polarization of the quadratic space $(B,\p)$ and thus
the space $S:=\wedge^*\p_-$ becomes a module under the Clifford algebra
$Cl(B,\p)$.
Since $\CG$ is simply connected the homomorphism $K\ra {\rm SO}(\p)$ given by
the adjoint action factors over the spin group ${\rm Spin}(\p)\subset
Cl(B,\p)$.
So $K$ acts on $S$.
The same applies to $\wedge^*\p_+$ and the $K$-action on the space
$\wedge^*\p_- \otimes \wedge^*\p_+ = \wedge^*\p$ coincides with the adjoint
action.

The action of $K$ on $S$ leaves invariant the subspaces
$$
S^+ := \wedge^{even} \p_- ,\ \ \ S^-:= \wedge^{odd}\p_-.
$$

Let $(\tau, V_\tau)$ be an irreducible unitary representation of $K$ then the
$K$-representations $S^\pm \otimes \tau$ define homogeneous vector bundles
$E(S^\pm \otimes \tau)$ over $X$ whose smooth sections can by identified with
the $K$-invariants:
$$
(C^\infty(G)\otimes S^\pm \otimes \tau)^K,
$$
where $K$ acts on $C^\infty(G)$ by right shifts.
The same applies to $S$ and we have $E(S\otimes\tau)=E(S^-\otimes\tau)\oplus
E(S^+\otimes\tau)$.
Let the Lie algebra $\g$ act on $C^\infty(G)$ by left invariant vector fields,
i.e. $Xf(g) := \frac{d}{dt}|_{t=0}f(g\exp(tX))$ for $X\in \g$ and $f\in
C^\infty(G)$.
Let $(X_j)_{1\leq j\leq \dim X}$ denote an orthonormal basis of $\p$ the we
have the {\bf Dirac operator} acting on $C^\infty(E(S\otimes\tau))$:
$$
D_\tau:= \sum_{j=0}^{\dim X} X_j \otimes c(X_j) \otimes 1,
$$
where we have written $c(X_j)$ for the Clifford action of $X_j\in\p\subset
Cl(B,\p)$ on $S$.
Clearly $D_\tau$ commutes $C^\infty(E(S^+\otimes \tau))$ and
$C^\infty(E(S^-\otimes \tau))$ and we will write $D_\tau^\pm$ for the
restriction of $D_\tau$ to $C^\infty(E(S^\pm\otimes \tau))$.
Then $D_\tau^+$ and $D_\tau^-$ are adjoints of each other.

The homogeneous bundle $E(S\otimes \tau)$ pushes down to a bundle $E_\Ga
(S\otimes \tau)$ over $X_\Ga$ whose space of smooth sections can be identified
with $(C^\infty (\Ga \bs G)\otimes S\otimes \tau)^K$.

In \cite{Lab-Pseu} it is proven that there is a compactly supported smooth
function $g_\tau$ such that
$$
\tr\pi (g_\tau) = \dim(V_\pi \otimes S^+ \otimes V_{\breve{\tau}} )^K -
\dim(V_\pi \otimes S^- \otimes V_{\breve{\tau}} )^K.
$$

We want to show that $\tr\pi(g_\tau)$ vanishes for a principal series
representation $\pi$.
To this end let $P=MAN$ be a parabolic subgroup with $A\subset \exp(\p)$.
Let $(\xi ,V_\xi)$ denote an irreducible unitary representation of $M$ and
$e^\nu$ a quasicharacter of $A$.
Let $\pi_{\xi ,\nu}:= {\rm Ind}_P^G \xi \otimes e^{\nu +\rho_P}\otimes 1$.

\begin{lemma} \label{pivonggleichnull}
We have $\tr\pi_{\xi ,\nu}(g_\tau) =0$.
\end{lemma}

\prf
By Frobenius reciprocity we have
$$
\Hom_K(\ga ,\pi_{\xi ,\nu}|_K) \cong \Hom_{K_M}(\ga |_{K_M},\xi ),
$$
where $K_M := K\cap M$ the claim will follow from $S^+|_{K_M} \cong
S^-|_{K_M}$.
To prove this let $0\neq \omega \in \p_-$ be in the image of the projection of
$\a :=\Lie A$ to $\p_-$.
Then $K_M$ acts trivially on $\C \omega\subset \p_-$.
Let $W\subset\p_-$ be a $K_M$-complement to $\C\omega$ then
$\bigwedge^*\p_- = \bigwedge^*W \oplus \omega \wedge \bigwedge^*W$ and so
$S^+ =$ $\bigwedge^{even}W \oplus \omega \wedge \bigwedge^{odd}W$ $ \cong_{K_M}
\bigwedge^{odd}W \oplus \omega \wedge \bigwedge^{even}W$ $=S^-$.
\qed

For any unitary representation $\pi$ of $G$ we define the Dirac operator
$$
D_{\tau ,\pi} := \sum_{j=0}^{\dim X} \pi(X_j) \otimes c(X_j) \otimes 1
$$
acting on $(\pi^\infty \otimes S \otimes \tau)^K$.
Write $D_{\tau ,\pi}^\pm$ for the restriction to $(\pi^\infty \otimes S^\pm
\otimes \tau)^K$.
If $\ker D_{\tau ,\pi}^+$ and $\ker D_{\tau ,\pi}^-$ are finite dimensional we
define
$$
\ind D_{\tau ,\pi}^+ := \dim \ker D_{\tau ,\pi}^+ -\dim \ker D_{\tau ,\pi}^-.
$$

\begin{lemma}
For $\pi \in \hat{G}$ the kernel of $D_{\tau ,\pi}$ is finite dimensional and
we have
$$
\ind D_{\tau ,\pi}^+ = \tr\pi (g_{\breve{\tau}}).
$$
\end{lemma}

\prf
Since $K$-types have finite multiplicities in $\pi$ it follows that
$(\pi^\infty \otimes S\otimes \tau)^K$ is finite dimensional.
On this finite dimensional space the operators $D_{\tau ,\pi}^+$ and $D_{\tau
,\pi}^-$ are adjoints of each other, so
\begin{eqnarray*}
\ind D_{\tau ,\pi}^+ &=& \dim \ker D_{\tau ,\pi}^- D_{\tau ,\pi}^+ - \dim \ker
D_{\tau ,\pi}^+ D_{\tau ,\pi}^-\\
	&=& \dim \ker D_{\tau ,\pi}^2 |_{(\pi^\infty \otimes S^+\otimes \tau)} - \dim
\ker D_{\tau ,\pi}^2 |_{(\pi^\infty \otimes S^-\otimes \tau)}.
\end{eqnarray*}

The formula of Parthasarathy \cite{Part} (see also \cite{AtSch}) implies
$$
D_{\tau ,\pi}^2 = -\pi(C) +\tau(C_K)+B(\rho_K)-B(\rho)|_{(\pi^\infty\otimes
S\otimes \tau)^K},
$$
where $C$ and $C_K$ are the Casimir operators of $G$ and $K$.

This gives the claim.
\qed

The Dirac operator for the $G$-representation on $L^2(\Ga \bs G)$ will be
denoted $D_{\tau ,\Ga}$.
Recall that $L^2(\Ga \bs G)$ decomposes as
$$
L^2(\Ga \bs G) = L^2(\Ga \bs G)_{disc} \oplus L^2(\Ga \bs G)_{cont},
$$
where $L^2(\Ga \bs G)_{disc}$, the discrete part, is the sum of all irreducible
subrepresentations of $L^2(\Ga \bs G)$ and $L^2(\Ga \bs G)_{cont}$ is a
continuous Hilbert integral
extended over the principal series.
Let $R(g_\tau)$ denote the convolution operator $\ph \mapsto \ph *
\breve{g}_\tau$, where $\breve{g}_\tau (x) :=g_\tau(x^{-1})$.
Then $R(g_\tau) = R_{cont}(g_\tau)+R_{disc}(g_\tau)$ accordingly and Lemma
\ref{pivonggleichnull} implies that $R_{cont}(g_\tau)=0$.
Further Theorem 7.1 of \cite{Art-inv2} says that $R(g_\tau)=R_{disc}(g_\tau)$
is a trace class operator. Therefore it follows

\begin{lemma} \label{index_is_trace_gamma}
The spaces $\ker D_{\tau ,\Ga}^+ \cap L^2(\Ga \bs G)$ and $\ker D_{\tau ,\Ga}^-
\cap L^2(\Ga \bs G)$ are finite dimensional.
Denote the difference of their dimensions by $\ind D_{\tau ,\Ga}^+$ then
$$
\ind D_{\tau ,\Ga}^+ = \tr (g_{\breve{\tau}} |L^2(\Ga \bs G)).
$$\qed
\end{lemma}

The association $\tau \mapsto g_\tau$ extends to virtual representations by
linearity.
Consider the virtual representation of $K$ on $S^+-S^-$.
We define
$$
f_\tau := g_{\tau \otimes (S^+ -S^-)}.
$$
It follows that for $\pi\in\hat{G}$ we have
$$
\tr\pi f_\tau = \sum_{q=0}^{\dim X} (-1)^q \dim (V_\pi \otimes \wedge^q \p
\otimes V_{\breve{\tau}})^K.
$$

An element $g$ of $G$ is called {\bf elliptic} if it lies in a compact subgroup
of $G$.
For any $g\in G$ and a compactly supported smooth function $f$ on $G$ let
$$
\CO_g(f) := \int_{G/G_g} f(xgx^{-1}) dx
$$
denote the {\bf orbital integral}.
The required normalization of Haar measures of $G$ and $G_g$ follows
\cite{HC-HA1}.

{}From \cite{Holtors} we take

\begin{proposition} \label{orbitalint}
Let $g$ be a semisimple element of the group $G$. If $g$ is not elliptic, the
orbital integrals $\CO_g(f_\sigma)$ and $\CO_g(g_\tau)$ vanish. If $g$ is
elliptic we may assume $g\in T$, where $T$ is a Cartan in $K$ and then we have
$$
\CO_g(f_\sigma) = \frac{{\tr\ \sigma(g)}|W(\t ,\g_g)| \prod_{\alpha \in
\Phi_g^+}(\rho_g ,\alpha)}{[G_g:G_g^0]c_g},
$$
for all elliptic $g$ and
$$
\CO_g(g_\tau) = \frac{{\tr\ \tau(g)}}{\det(1-g^{-1} | \p_+)},
$$
if $g$ is regular elliptic. For general elliptic $g$ we have
$$
\CO_g(g_\tau) = \frac{\sum_{s\in W(T,K)} \det(s) \tilde{\omega}_g
g^{s\la_{{\tau}}+\rho -\rho_K}}
        {[G_g :G_g^0]c_g g^\rho \prod_{\alpha \in \phi^+ -
\phi_g^+}(1-g^{-\alpha})},
$$
where $c_g$ is Harish-Chandra's constant, it does only depend on the
centralizer $G_g$ of g. Its value is given in \cite{D-Hitors}, further
$\tilde{\omega}_\ga$ is the differential operator as in \cite{HC-DS} p.33.
\qed
\end{proposition}

\begin{corollary} \label{orbintvanish}
If $G$ is a product of real rank one groups then the orbital integral
$\CO_g(f_\tau)$ vanishes also for $g$ non-semisimple.
\end{corollary}

\prf
It suffices to assume that the real rank of $G$ is one.
Consider a non-semisimple element $g$. By \cite{Barb}, sec. 6 we get a curve
$t\mapsto z_t$ of semisimple elements and a natural number $m$ such that
$\CO_g(f_\tau)=\lim_{t\rightarrow 0} t^{m/2}\CO_{z_t}(f_\tau)$.
The function $x\mapsto \CO_x(f_\tau)$ is bounded on semisimple elements by the
proposition and therefore we have $\CO_g(f_\tau) =0$ for $g$ non semisimple.
\qed

\begin{conjecture} \label{all_orbitalint_vanish}
The above corollary holds for all groups $G$.
\qed
\end{conjecture}

\section{The $G$-index}
Let $\CR_G \subset B(L^2(G))$ denote the von Neumann algebra defined as the
commutant of the left representation of $G$ on $L^2(G)$.
On $\CR_G$ there is a canonical faithful, normal, semi-finite trace $\tr_G$,
called the {\bf $G$-trace}, uniquely determined by the property that
$$
\tr_G(R(f)^*R(f))=\int_G |f(g)|^2 dg,
$$
where $R$ denotes the right representation of $G$
(See \cite{ConMosc}).
Let $\dim_G$ denote the dimension defined by the $G$-trace.

Choose a compact form $G^d$ of $G$ in a way that $G^d$ contains the compact
group $K$.
The homogeneous space $X^d = G^d/K$ then is symmetric and is called the {\bf
dual symmetric space} to $X$.

Let $\la =\la_\tau$ denote the infinitesimal character of $\tau$ then $\la$
also defines an infinitesimal character of some irreducible representation
$W_\la$ of $G^d$.

\begin{proposition}
The spaces $\ker D_\tau^\pm \cap L^2(E(S^\pm\otimes \tau))$ are
finite-dimensional under $\dim_G$.
Let {\bf $\ind_G(D_\tau^+)$} denote the {\bf $G$-index} of $D_\tau^+$, which
is, by definition, the difference of these two $G$-dimensions.
Then if $\la$ is regular with respect to the full root system we have
$$
\ind_G(D_\tau^+) = \frac{\dim W_\la}{\chi(X^d)},
$$
where $\chi(X^d)$ is the Euler-characteristic of $X^d$.
If $\la$ is not regular, the $G$-index of $D_\tau^+$ vanishes.
\end{proposition}

It is known that the Euler characteristic is positive and that
$\chi(X^d)=|W(T,G^d)|$, where $T$ is a Cartan subgroup of $G^d$.

\prf
The finite dimensionality is Lemma 3.2 in \cite{ConMosc}.
The index formula follows from formulas (3.7)-(3.13) of \cite{AtSch}.
\qed

\section{The index theorem}

Assume now that $\CG$ is the restriction to $\Q$ of some algebraic group $\CH$
over a number field $F$.
It is known that the index of the Dirac operator $D_{\tau ,\Ga}$ vanishes if
$\rank G>\rank K$.
If $F$ has a complex place then it follows $\rank G>\rank K$.
So, in order to have a nontrivial theory we will assume the field $F$ to be
totally real.
Let $S$ denote the set of archimedian places of $F$.
We assume that $S$ has at least two elements.

\begin{theorem}
The index of the operator $D_{\tau ,\Ga}^+$ is given by
$$
\ind (D_{\tau ,\Ga}^+) = \ind_G (D_\tau^+) \chi(X_\Ga) + R(\tau ,\Ga),
$$
where the "error term" $R(\tau ,\Ga)$ equals
$$
R(\tau ,\Ga)= \sum_{\ga \in \Ga_{ns}/G} a(\ga) \CO_\ga (g_\tau).
$$
The sum is extended over the set of non-semisimple elements $\Ga_{ns}$ of $\Ga$
modulo the equivalence relation $(G,S)$ defined in \cite{Art-inv2}
The constant $a(\ga)$ coincides with $a^G(S,\ga)$ of \cite{Art-inv2} up to a
volume factor.
\end{theorem}

Note that the first summand on the right hand side of the index formula also
coincides with the integral $\int_{X_\Ga}\omega$, where $\omega$ is the index
form of $D_\tau$. This also equals the $\Ga$-index \cite{atiyah} of $D_\tau^+$,
so we have
$\ind(D_{\tau ,\Ga}^+)=\ind_\Ga (D_\tau^+) + R(\tau ,\Ga)$.

\prf
Let $K_\Ga$ be the compact open subgroup of $\CG(\A_{fin})$ such that $\Ga
=K_\Ga \cap \CG(\Q)$.
Define a compactly supported function $f$ on $\CG(\A)$ by $f=f_{fin}\otimes
f_\infty$, where $f_{fin}:=\frac{1}{\vol(K_\Ga)}{\bf 1}_{K_\Ga}$ and $f_\infty
:=g_\tau$.
Plug the function $f$ into  Theorem 7.1 (b) of \cite{Art-inv2}.
Then use Corollary \ref{orbintvanish}.
\qed

Unfortunately the constants $a(\ga)$ are not easy to compute (see the comment
on page 209 of \cite{Art-orbits}).
In simple cases like $\CG = \Res_{F/\Q}SL_2$ however,one will probably be able
to get more explicit expressions.
In this paper we will be more interested in the vanishing of $R(\tau ,\Ga)$.

\begin{corollary}
In the case $G= \Res_{F/G}SO(n,1)$ for $n\geq 3$, the error term vanishes.
Here $SO(n,1)$ stands for the special orthogonal group of a quadratic form
which has signature $(n,1)$ over the reals.
\end{corollary}

This extends Thm 7.5 of \cite{BM}.

\prf
Proposition \ref{orbitalint} and sec. 6 of \cite{Barb} imply that $R(\tau
,\Ga)=0$ if $K$ has discrete center.
The latter condition is satisfied for $SO(n,1)$.
\qed

\section{The Euler operator}
Consider the homogeneous vector bundle $E=E(\tau)$ associated with the
representation $\tau$.
Choose a homogeneous connection on $E(\tau)$.
Those always exist and in case that $\tau=\sigma |_K$, where $\sigma$ is a
finite dimensional representation of $G$ there is a unique flat homogeneous
connection.
This is the case considered in \cite{stern2}.
The choice of a connection gives us an exterior differential
$$
d : \Omega^.(E) \ra \Omega^{.+1}(E).
$$
The representation $\tau$ being unitary gives us a homogeneous hermitian metric
on $E$ and we can define the formal adjoint $d^*$ of $d$ and the {\bf Euler
operator}:
$$
D_\tau^e := d+d^* : \Omega^{even}(E) \ra \Omega^{odd}(E).
$$

We now come to the main result of this paper:

\begin{theorem}
Assume that $G=\Res_{F/\Q}\CH$, where $\CH$ is a semisimple connected linear
algebraic group over the totally real number field $F$.
Assume $F$ has at least two real embeddings and $\CH$ is over $\R$ a product of
rank one groups.
Let $G:=\CG(\R)$ and
$X=G/K$ the symmetric space attached to $G$.
Let $\Ga$  a torsion free congruence subgroup of $G$. Write $X_\Ga =\Ga \bs X$
for the quotient manifold.
Let $E(\tau)$ be a homogeneous vector bundle over $X$ given by a unitary finite
dimensional representation $\tau$ of the compact group $K$ and let $D_{\tau
,\Ga}^e$ be the Euler operator of the pushdown of $E(\tau)$ to $X_\Ga$ then
$D_{\tau ,\Ga}^e$ has a well defined index and
$$
\ind (D_{\tau ,\Ga}^e) = \ind_G(D_\tau^e) \chi(X_\Ga).
$$
\end{theorem}

\begin{corollary}
Assume Conjecture \ref{all_orbitalint_vanish} holds. Then the condition on $G$
in the above theorem can be removed.
\end{corollary}

\prf
The index of $D_{\tau ,\Ga}^e$ coincides with the index of $D_{\tilde{\tau}
,\Ga}^+$, where $\tilde{\tau}$ is the virtual representation $\tau \otimes (S^+
-S^+)$.
With Lemma \ref{index_is_trace_gamma} it follows that the index exists.
The index theorem tells us that
$$
\ind (D_{\tau ,\Ga}^e) = \ind_G(D_\tau^e) \chi(X_\Ga) + R(\tilde{\tau},\Ga).
$$
Corollary \ref{orbintvanish} tells us that $R(\tilde{\tau},\Ga)=0$.
\qed

\begin{corollary}
For the usual Euler operator $D_\Ga^e$ on $X_\Ga$ we have
$$
\ind (D_{\Ga}^e) = \chi(X_\Ga).
$$
\end{corollary}

\prf
This corollary is clearly valid for cocompact groups $\Ga$. This implies that
the $G$-index must be one.
\qed

\tiny
Version: \today

\end{document}